\documentclass{article}
\usepackage{graphicx}
\usepackage{amsmath}
\usepackage{amssymb}
\usepackage{amsthm}

\title{Radial Limits of Bounded Nonparametric PMC Surfaces}
\author{Mozhgan Entekhabi \& Kirk E. Lancaster
                                  \\
                       Department of Mathematics, Statistics \& Physics \\
                            Wichita State University \\
                            Wichita, Kansas, 67260-0033}

\def\Real{{\rm I\hspace{-0.2em}R}}

\newcommand\myeq{\mathrel{\overset{\makebox[0pt]{\mbox{\normalfont\tiny\sffamily def}}}{=}}}

\date{ }

\begin{document}
\maketitle

\vspace*{3mm}

\begin{center}
  \Large \emph{Dedicated to the memory of Alan Ross Elcrat}
\end{center}

\vspace*{3mm}

\begin{abstract}

\noindent Consider a solution $f\in C^{2}(\Omega)$  of a prescribed mean curvature equation 
\[
{\rm div}\left(\frac{\nabla f}{\sqrt{1+|\nabla f|^{2}}}\right)=2H(x,f) \ \ \ \ {\rm in} \ \ \Omega,
\]
where $\Omega\subset \Real^{2}$  is a domain whose boundary has a corner at ${\cal O}=(0,0)\in\partial\Omega.$  
If $\sup_{x\in\Omega} |f(x)|$  and $\sup_{x\in\Omega} |H(x,f(x))|$  are both finite and 
$\Omega$  has a reentrant corner at ${\cal O},$   then the radial limits of $f$  at ${\cal O},$  
\[
Rf(\theta) \myeq \lim_{r\downarrow 0} f(r\cos(\theta),r\sin(\theta)), 
\]
are shown to exist and to have a specific type of behavior, independent of the boundary behavior of $f$  on 
$\partial\Omega\setminus \{ {\cal O} \}.$  
If $\sup_{x\in\Omega} |f(x)|$  and $\sup_{x\in\Omega} |H(x,f(x))|$  are both finite and the trace of $f$  on one side has a 
limit at ${\cal O},$  then the radial limits of $f$  at ${\cal O}$  exist and have a specific type of behavior. 
\end{abstract}
\newtheorem{thm}{Theorem}
\newtheorem{prop}{Proposition}
\newtheorem{cor}{Corollary}
\newtheorem{lem}{Lemma}
 
\section{Introduction and Statement of Main Theorems}
 
\begin{figure}[ht]
\centering
\includegraphics{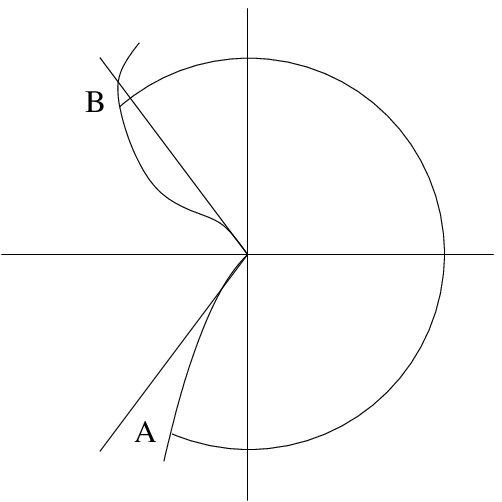}
\caption{The domain $\Omega^{*}$  \label{ZERO}}
\end{figure}

Consider the prescribed mean curvature equation
\begin{equation}
\label{PMC}
Nf  =  2H(\cdot,f)  \mbox{  \ in \ } \Omega    
\end{equation}
where $\Omega$ is a domain  in ${\Real}^{2}$  whose boundary has a corner at ${\cal O}\in\partial\Omega,$  
$Nf = \nabla \cdot Tf = {\rm div}\left(Tf\right),$  $Tf = \frac{\nabla f}{\sqrt{1 + |\nabla f|^{2}}},$   
$H:\Omega\times\Real \to \Real$  and $H$  satisfies one of the conditions which guarantees that ``cusp solutions'' 
(e.g. \S 5 of \cite{LS1}, \cite{LS2}) do not exist; for example,  $H({\bf x},t)$  is strictly increasing in $t$  
for each ${\bf x}$  or is real-analytic (e.g. constant).  We will assume ${\cal O}=(0,0).$  
Let $\Omega^{*} = \Omega \cap B_{\delta^{*}}({\cal O})$, where $B_{\delta^{*}}({\cal O})$ is the ball in $\Real^{2}$ 
of radius 
$\delta^{*}$ about ${\cal O}$.   Polar coordinates relative to ${\cal O}$ will be denoted by $r$ and $\theta$.
We assume that $\partial \Omega$ is piecewise smooth and there exists $\alpha \in \left(0,\pi\right)$  
such that $\partial \Omega \cap B_{\delta^{*}}({\cal O})$  consists of two arcs  ${\partial}^{+}\Omega^{*}$  
and $\partial^{-}\Omega^{*}$, whose tangent lines approach the lines $L^{+}: \:  \theta = \alpha$  and 
$L^{-}: \: \theta = - \alpha$, respectively, as the point ${\cal O}$ is approached (see Figure 1 of \cite{LS1errata} 
or Figure \ref{ZERO}).  

\noindent Suppose  
\begin{equation}
\label{Bounds}
\sup_{x\in\Omega} |f(x)| < \infty \ \ \ \  {\rm and} \ \ \ \  \sup_{x\in\Omega} |H(x,f(x))| < \infty.  
\end{equation}
We shall prove 

\begin{thm}
\label{CONCLUSION}
Let $f\in C^{2}(\Omega)$  satisfy (\ref{PMC}) and suppose (\ref{Bounds}) holds and $\alpha \in \left(\frac{\pi}{2},\pi\right).$
Then for each $\theta\in (-\alpha,\alpha),$ 
\[
Rf(\theta) \myeq \lim_{r\downarrow 0} f(r\cos(\theta),r\sin(\theta)) 
\]
exists and $Rf(\cdot)$  is a continuous function on $(-\alpha, \alpha)$ which behaves in one of the following ways:

\noindent (i)  $Rf:(-\alpha,\alpha)\to\Real$  is a constant function (i.e. $f$  has a nontangential limit at ${\cal O}$). 

\noindent (ii) There exist $\alpha_{1}$ and $\alpha_{2}$ so that $-\alpha \leq \alpha_{1}
< \alpha_{2} \leq \alpha$ and $Rf$ is constant on $(-\alpha, \alpha_{1}]$ and
$[ \alpha_{2}, \alpha)$ and strictly increasing or strictly decreasing on
$(\alpha_{1}, \alpha_{2})$.  
 
\noindent (iii) There exist $\alpha_{1}, \alpha_{L}, \alpha_{R}, \alpha_{2}$ so that
$-\alpha \leq \alpha_{1} < \alpha_{L} < \alpha_{R} < \alpha_{2} \leq \alpha,
\alpha_{R}= \alpha_{L} + \pi$, and $Rf$ is constant on $(-\alpha, \alpha_{1}],
[ \alpha_{L}, \alpha_{R}]$, and $[ \alpha_{2}, \alpha)$ and either  strictly increasing
on $(\alpha_{1}, \alpha_{L}]$ and  strictly decreasing on $[ \alpha_{R}, \alpha_{2})$ or
 strictly decreasing on $(\alpha_{1}, \alpha_{L}]$ and  strictly increasing on $[\alpha_{R},\alpha_{2})$.  
 \end{thm}

\noindent At a convex corner (i.e. $\alpha \in \left(0,\frac{\pi}{2}\right]),$   Theorem \ref{CONCLUSION} is not 
applicable.   The additional assumption that the trace of $f$  on one side (e.g. $\partial^{-}\Omega^{*}$) 
has a limit at ${\cal O}$  implies the radial limits of $f$  exist.

\begin{thm}
\label{DEMONSTRATION}
Let $f\in C^{2}(\Omega)\cap C^{0}\left(\Omega\cup \partial^{-}\Omega^{*}\setminus \{ {\cal O} \}\right)$  
satisfy (\ref{PMC}).  Suppose  (\ref{Bounds}) holds  and 
$m=\lim_{\partial^{-}\Omega^{*}\ni {\bf x}\to {\cal O} } f\left({\bf x}\right)$  exists.
Then for each $\theta\in (-\alpha,\alpha),$   $Rf(\theta)$  exists and $Rf(\cdot)$  is a continuous function 
on $[-\alpha, \alpha),$  where $Rf(-\alpha) \myeq m.$  
If  $\alpha \in \left(0,\frac{\pi}{2}\right],$  $Rf$  can behave as in (i) or (ii) in Theorem \ref{CONCLUSION}.  
If $\alpha \in \left(\frac{\pi}{2},\pi \right),$  $Rf$  can behave as in (i), (ii) or (iii) in Theorem \ref{CONCLUSION}. 
\end{thm}

\noindent The conclusions of these theorems were first obtained in \cite{Lan:85} for minimal surfaces satisfying Dirichlet 
boundary conditions and then for nonparametric prescribed mean curvature surfaces satisfying Dirichlet (\cite{EL:86,Lan:88}) 
or contact angle (\cite{LS1}) boundary conditions; see also \cite{JL1,Lan:91}.  
Notice that Theorem \ref{CONCLUSION} applies to a solution of a capillary surface problem whose domain has a reentrant corner 
even when the contact angle equals $0$  and/or $\pi$  on some (or all) of $\partial  \Omega^{*}.$  
\vspace{3mm}

\noindent {\bf Remark:} {\it Notice that the assumption that $\Omega$  has a reentrant corner at ${\cal O}\in\partial\Omega$  
or that the trace of $f$  from one side of $\partial\Omega$  is continuous at ${\cal O}$  is critical here; 
the nonexistence of radial limits at $(1,0)$  when $\Omega=B_{1}\left( {\cal O} \right)$  and the boundary data is symmetric 
with respect to the horizontal axis is demonistrated in \cite{Lan:89} and in Theorem 3 of \cite{LS1}.  
In \cite{Lan:87}, the second author conjectured the existence of radial limits at corners for 
bounded solutions of Dirichlet problems for the minimal surface equation in $\Real^{2},$  independent of boundary conditions.  
Although \cite{Lan:89} proved this conjecture false, Theorems \ref{CONCLUSION} and \ref{DEMONSTRATION} show it is true 
in many cases.}

\section{Proof of Theorem \ref{CONCLUSION}}  

Since $f\in C^{2}(\Omega)$  (and so in $C^{0}(\Omega)$), we may assume that $f$  is uniformly continuous on 
$\{{\bf x}\in \Omega^{*} \ :\ |{\bf x}|>\delta\}$  for each $\delta\in (0,\delta^{*});$  
if this is not true, we may replace $\Omega$  with $U,$  $U\subset \Omega,$  such that 
$\partial\Omega \cap \partial U = \{ {\cal O} \}$  and $\partial U \cap B_{\delta^{*}}({\cal O})$  
consists of two arcs  ${\partial}^{+} U$  and $\partial^{-} U$, whose tangent lines approach the lines 
$L^{+}: \:  \theta = \alpha$  and $L^{-}: \: \theta = - \alpha$, respectively, as the point ${\cal O}$ is approached.
Set
\[ 
S^{*}_{0} = \{ ({\bf x},f({\bf x})) : {\bf x} \in \Omega^{*} \}
\]
and
\[ 
\Gamma^{*}_{0} = \{ ({\bf x},f({\bf x})): {\bf x} \in \partial \Omega^{*}    \setminus \{ {\cal O} \} \}; 
\]
the points where $\partial B_{\delta^{*}}({\cal O})$ intersect $\partial \Omega$  are labeled $A\in {\partial}^{-}\Omega^{*}$ and 
$B \in {\partial}^{+}\Omega^{*}.$
From the calculation on page 170 of \cite{LS1},  we see that the area of $S^{*}_{0}$  is finite; let $M_{0}$  denote this area. 
For $\delta\in (0,1),$  set 
\[
p(\delta) = \sqrt{\frac{8\pi M_{0}}{\ln\left(\frac{1}{\delta}\right)}}.
\]
Let $E= \{ (u,v) : u^{2}+v^{2}<1 \}.$ 
As in \cite{EL:86,LS1}, there is a parametric description of the surface $S^{*}_{0},$  
\begin{equation}
\label{PARAMETRIC}
Y(u,v) = (a(u,v),b(u,v),c(u,v)) \in C^{2}(E:{\Real}^{3}), 
\end{equation}
which has the following properties:
 
\noindent $\left(a_{1}\right)$  $Y$ is a diffeomorphism of $E$ onto $S^{*}_{0}$.
 
\noindent $\left(a_{2}\right)$
Set $G(u,v)=(a(u,v),b(u,v)),$  $(u,v)\in E.$   Then $G \in C^{0}(\overline{E} : {\Real}^{2}).$  

\noindent $\left(a_{3}\right)$
Let $\sigma=G^{-1}\left(\partial \Omega_{1}\setminus \{ {\cal O} \}\right);$  
then $\sigma$ is a connected arc of $\partial E$  and $Y$ maps $\sigma$ strictly monotonically onto 
$\Gamma^{*}_{0}.$  
We may assume the endpoints of $\sigma$  are ${\bf o}_{1}$  and ${\bf o}_{2}$  and there exist points 
${\bf a}, {\bf b}\in\sigma$  such that $G({\bf a})=A,$  $G({\bf b})=B,$
$G$  maps the (open) arc ${\bf o}_{1}{\bf b}$  onto $\partial^{+}\Omega,$  and $G$ 
maps the (open) arc ${\bf o}_{2}{\bf a}$  onto $\partial^{-}\Omega.$ 
(Note that ${\bf o}_{1}$  and ${\bf o}_{2}$  are not assumed to be distinct at this point; 
one of Figure 4a or 4b of \cite{LS1errata} illustrates this situation.)   
 
\noindent $\left(a_{4}\right)$
$Y$ is conformal on $E$: $Y_{u} \cdot Y_{v} = 0, Y_{u}\cdot Y_{u} = Y_{v}\cdot Y_{v}$
on $E$.
 
\noindent $\left(a_{5}\right)$
$\triangle Y := Y_{uu} + Y_{vv} = H\left(Y\right)Y_{u} \times Y_{v}$  on $E$.
\vspace{1mm} 

\noindent Here by the (open) arcs ${\bf o}_{1}{\bf b}$  and ${\bf o}_{2}{\bf a}$  are meant the 
component of $\partial E\setminus\{{\bf o}_{1},{\bf b}\}$  which does not contain  ${\bf a}$ 
and the component of $\partial E\setminus\{{\bf o}_{2},{\bf a}\}$  which does not contain  ${\bf b}$  
respectively. 
Let  $\sigma_{0} = \partial E \setminus \sigma$.
\vspace{2mm}

There are two cases we wish to consider: 
\begin{itemize}
\item[$\left(A\right)$] ${\bf o}_{1}= {\bf o}_{2}.$
\item[$\left(B\right)$] ${\bf o}_{1}\neq {\bf o}_{2}.$
\end{itemize}
These correspond to Cases 5 and 3 respectively in Step 1 of the proof of Theorem~1 of \cite{LS1}. 

\begin{figure}[ht]
\centering
\includegraphics{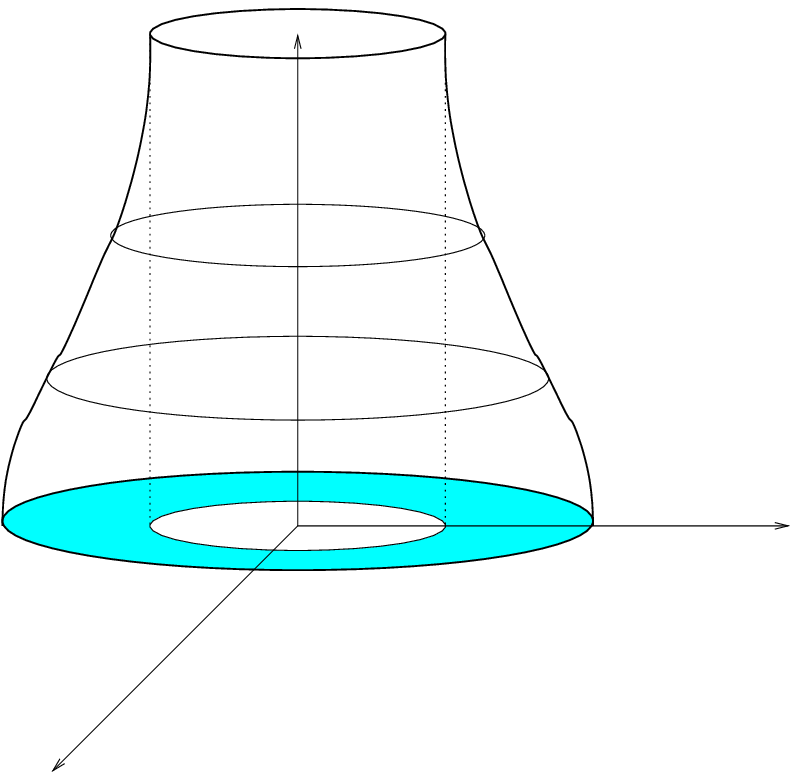}
\caption{The graph of $h$  over ${\cal A}$  \label{FOUR}}
\end{figure}

Let us first assume that $\left(A\right)$  holds  and set ${\bf o}={\bf o}_{1}={\bf o}_{2}.$
Let $h$ denote a function on  the annulus ${\cal A}=\{ {\bf x} : r_{1} \leq |{\bf x}| \leq r_{2} \}$  
which vanishes on the circle $|{\bf x}|=r_{2}$  and whose graph is an unduloid surface  with constant mean 
curvature $-H_{0}$  which  becomes vertical at $|{\bf x}| = r_{1}$  and at $|{\bf x}| = r_{2}$  
(see Figure \ref{FOUR})  for suitable 
$r_{1}<r_{2}$  (e.g. \cite{LS1}, pp. 170-1).  Let $q$  denote the modulus of continuity 
of $h$  (i.e. $|h({\bf x}_{1}) - h({\bf x}_{2})| \leq q(| {\bf x}_{1} - {\bf  x}_{2}|).$
For each ${\bf p}\in \Real^{2}$  with $|{\bf p}|=r_{1},$  set 
${\cal A}({\bf p})=\{ {\bf x} : r_{1} \leq |{\bf x}-{\bf p}|\leq r_{2} \}$   and define 
$h_{{\bf p}}:{\cal A}({\bf p})\to\Real$  by $h_{{\bf p}}({\bf x})=h\left({\bf x}-{\bf p}\right).$

For $r>0,$  set $B_{r}= \{ {\bf u} \in \overline{E} : |{\bf u} - {\bf o}| < r \},$   
$C_{r} = \{ {\bf u} \in \overline{E} : |{\bf u} - {\bf o}| = r \}$  and let $l_{r}$ be the length of the image curve 
$Y(C_{r});$  also let $C^{\prime}_{r} = G(C_{r})$  and $B^{\prime}_{r}= G(B_{r}).$   
From the Courant-Lebesgue Lemma (e.g. Lemma $3.1$ in \cite{Cour:50}), we see that 
for each $\delta\in (0,1),$ 
there exists a $\rho=\rho(\delta)\in \left(\delta,\sqrt{\delta}\right)$  such that the arclength 
$l_{\rho}$  of $Y(C_{\rho})$  is less than $p(\delta).$   
For $\delta>0,$  let $k(\delta)= \inf_{{\bf u}\in C_{\rho(\delta)}}c({\bf u}) = \inf_{ {\bf x}\in C^{\prime}_{\rho(\delta)} } f({\bf x})$  
and                  $m(\delta)= \sup_{{\bf u}\in C_{\rho(\delta)}}c({\bf u}) = \sup_{ {\bf x}\in C^{\prime}_{\rho(\delta)} } f({\bf x});$ 
notice that $m(\delta)-k(\delta)\le l_{\rho} < p(\delta).$  

\begin{figure}[ht]
\centering
\includegraphics{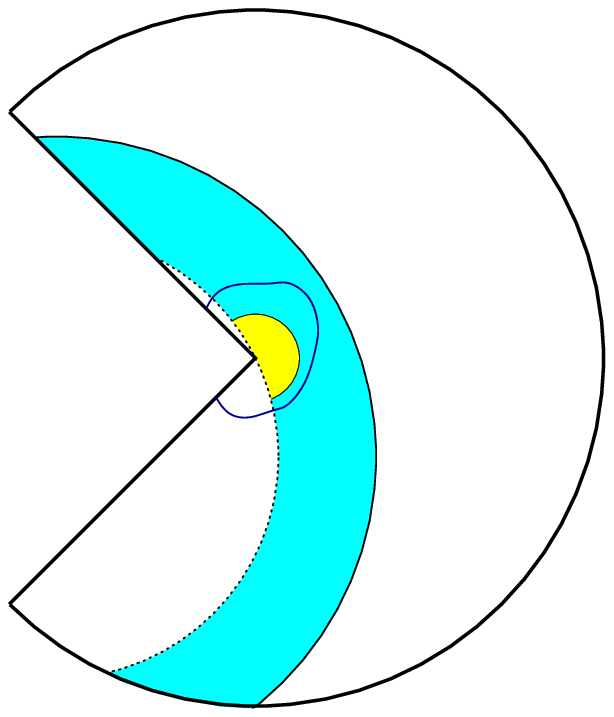}
\caption{$\Omega^{*}\cap {\cal A}\left({\bf p}_{1}\right),$  $C^{\prime}_{\rho(\delta)}$  (blue curve),   
$B_{\eta(\delta)}({\cal O})\cap{\cal A}\left({\bf p}_{1}\right)$ (yellow) \label{THREE}}
\end{figure}

For each $\delta\in (0,1)$  with $\sqrt{\delta}<\min\{|{\bf o}-{\bf a}|, |{\bf o}-{\bf b}|\},$  there are two points in 
$C_{\rho(\delta)}\cap\partial E;$   we denote these points as ${\bf e}_{1}(\delta)\in {\bf o}{\bf b}$  and 
${\bf e}_{2}(\delta)\in {\bf o}{\bf a}$  and set ${\bf y}_{1}(\delta)=G({\bf e}_{1}(\delta))$  
and ${\bf y}_{2}(\delta)=G({\bf e}_{2}(\delta)).$
Notice that $C^{\prime}_{\rho(\delta)}$  is a curve in $\overline{\Omega}$  which joins ${\bf y}_{1}\in {\partial}^{+}\Omega^{*}$ 
and ${\bf y}_{2}\in {\partial}^{-}\Omega^{*}$  and 
$\partial\Omega\cap C^{\prime}_{\rho(\delta)}\setminus \{{\bf y}_{1},{\bf y}_{2}\}=\emptyset;$  
therefore there exists $\eta=\eta(\delta)>0$  such that  
$B_{\eta(\delta)}({\cal O})=\{ {\bf x}\in \Omega : |{\bf x}|<\eta(\delta)\} \subset B^{\prime}_{\rho(\delta)}$  
(see Figure \ref{THREE}).

\begin{figure}[ht]
\centering
\includegraphics{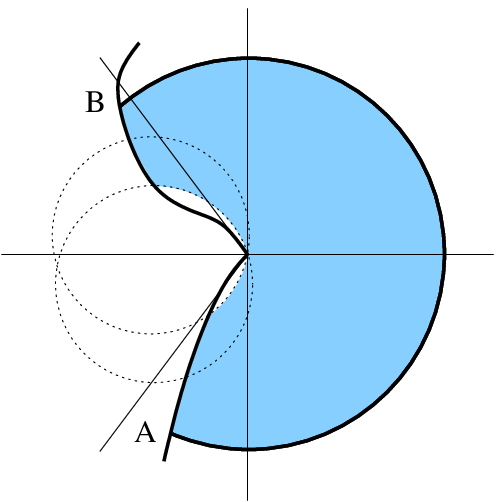}
\caption{$\Omega_{0}$ \label{TWO}}
%\caption{$\Omega^{*}\cap {\cal A}({\bf p}_{1})$ \label{ONE}}
\end{figure}

Fix $\delta_{0}\in (0,\delta^{*})$  with $\sqrt{\delta_{0}}<\min\{|{\bf o}-{\bf a}|, |{\bf o}-{\bf b}|\}.$  
Let ${\bf p}_{1}\in\Real^{2}$  satisfy $|{\bf p}_{1}|=r_{1}$  and $|{\bf p}_{1}-{\bf y}_{1}(\delta_{0})|=r_{1}$  such that ${\bf p}_{1}$  
lies below (and to the left of) the line through ${\cal O}$  and ${\bf y}_{1}(\delta_{0}).$  
Let ${\bf p}_{2}\in\Real^{2}$  satisfy $|{\bf p}_{2}|=r_{1}$  and $|{\bf p}_{2}-{\bf y}_{2}(\delta_{0})|=r_{1}$  such that ${\bf p}_{2}$  
lies above (and to the left of) the line through ${\cal O}$  and ${\bf y}_{2}(\delta_{0}).$  
Set $\Omega_{0} = \{ {\bf x}\in \Omega^{*} \ : \ |{\bf x}-{\bf p}_{1}|>r_{1} \} \cup 
\{ {\bf x}\in \Omega^{*} \ : \ |{\bf x}-{\bf p}_{2}|>r_{1} \}$  
(see Figure \ref{TWO}).
\vspace{1mm}

\noindent {\bf Claim:} $f$  is uniformly continuous on $\Omega_{0}.$
\vspace{1mm}

\noindent {\bf Pf:}  Let $\epsilon>0.$  Choose $\delta\in (0,\delta_{0})$  such that  $p(\delta)+q(p(\delta))<\frac{\epsilon}{4}$  
and $p(\delta)<r_{2}-r_{1}.$    
Pick a point ${\bf w}\in C^{\prime}_{\rho(\delta)}$  and define $b^{\pm}_{j}:{\cal A}({\bf p}_{j})\to\Real$  by 
\[
b^{\pm}_{j}({\bf x})=f({\bf w})  \pm p(\delta)  \pm h_{{\bf p}_{j}}({\bf x}), \ \ \ \ {\bf x}\in {\cal A}({\bf p}_{j})
\]
for $j\in \{1,2\}.$  Notice that 
\[
b^{-}_{j}({\bf x})<f({\bf x})<b^{+}_{j}({\bf x}) \ \ \ \ {\rm for} \ \  
{\bf x}\in B^{\prime}_{\rho(\delta)}  \cap {\cal A}({\bf p}_{j}), \ \ \ \ j\in \{1,2\}.
\]

\noindent If ${\bf x}_{1},{\bf x}_{2}\in \Omega_{0}$  satisfy  $|{\bf x}_{1}|<\eta(\delta)$  and $|{\bf x}_{2}|<\eta(\delta),$  
then there exist ${\bf x}_{3} \in {\cal A}({\bf p}_{1}) \cap {\cal A}({\bf p}_{2})$   with  $|{\bf x}_{3}|<\eta(\delta)$
such that $|{\bf x}_{1}-{\bf x}_{3}|<\eta(\delta)$  and $|{\bf x}_{2}-{\bf x}_{3}|<\eta(\delta)$  and so 
\begin{equation}
\label{Peach}
|f({\bf x}_{1})-f({\bf x}_{2})| \le |f({\bf x}_{1})-f({\bf x}_{3})| + |f({\bf x}_{1})-f({\bf x}_{3})|
<4p(\delta)+4q\left(p(\delta)\right)<\epsilon.
\end{equation}
Since $f$  is uniformly continuous on $\{ {\bf x}\in \Omega^{*} \ : \ |{\bf x}|\ge \frac{1}{2}\eta(\delta)\},$  
there exists a $\lambda>0$  such that if ${\bf x}_{1},{\bf x}_{2}\in \Omega^{*}$  satisfy 
$|{\bf x}_{1}-{\bf x}_{2}|\ge \frac{1}{2}\eta(\delta)$  and 
$|{\bf x}_{1}-{\bf x}_{2}| < \lambda,$  then  $|f({\bf x}_{1})-f({\bf x}_{2})|<\epsilon.$
Now set $d=d(\epsilon)=\min\{\lambda, \frac{1}{2}\eta(\delta)\}.$  
If ${\bf x}_{1},{\bf x}_{2}\in \Omega_{0},$  $|{\bf x}_{1}-{\bf x}_{2}|< d(\epsilon)\le \frac{1}{2}\eta(\delta)$  
and $|{\bf x}_{1}|<\frac{1}{2}\eta(\delta),$  then $|{\bf x}_{1}|<\eta(\delta)$  and $|{\bf x}_{2}|<\eta(\delta);$  
hence $|f({\bf x}_{1})-f({\bf x}_{2})|<\epsilon$  by (\ref{Peach}).    
Next, if ${\bf x}_{1},{\bf x}_{2}\in \Omega_{0},$  $|{\bf x}_{1}-{\bf x}_{2}|< d(\epsilon)\le\lambda,$ 
$|{\bf x}_{1}|\ge \frac{1}{2}\eta(\delta)$  and $|{\bf x}_{2}|\ge \frac{1}{2}\eta(\delta),$
then $|f({\bf x}_{1})-f({\bf x}_{2})|<\epsilon.$  
Therefore, for all ${\bf x}_{1},{\bf x}_{2}\in \Omega_{0}$  with $|{\bf x}_{1}-{\bf x}_{2}|< d(\epsilon),$  
we have $|f({\bf x}_{1})-f({\bf x}_{2})|<\epsilon.$   The claim is proven.
\vspace{2mm}

If $\{\left(r\cos(\theta^{-}(\delta_{0})),r\sin(\theta^{-}(\delta_{0}))\right) : r\ge 0\}$  is the tangent ray to 
$\partial {\cal A}({\bf p}_{2})$  at ${\cal O},$  
$\{\left(r\cos(\theta^{+}(\delta_{0})),r\sin(\theta^{+}(\delta_{0}))\right) : r\ge 0\}$  is the tangent ray to 
$\partial {\cal A}({\bf p}_{1})$  at ${\cal O}$  and  
$\theta^{-}(\delta_{0}),\theta^{+}(\delta_{0})\in (-\alpha,\alpha),$  then it follows from the Claim that 
$f\in C^{0}(\overline{\Omega_{0}}),$  the radial limits $Rf(\theta)$  of $f$  at ${\cal O}$  exist for 
$\theta\in [\theta^{-}(\delta_{0}),\theta^{+}(\delta_{0})]$  and the radial limits are identical 
(i.e. $Rf(\theta)=f({\cal O})$  for all $\theta\in [\theta^{-}(\delta_{0}),\theta^{+}(\delta_{0})].$)
Since  
\begin{equation}
\label{Tree}
\lim_{\delta_{0}\downarrow 0} \theta^{-}(\delta_{0}) = -\alpha 
\ \ \ \   {\rm and} \ \ \ \ 
\lim_{\delta_{0}\downarrow 0} \theta^{+}(\delta_{0}) = \alpha,
\end{equation}
Theorem \ref{CONCLUSION} is proven in this case.
\vspace{2mm}

Let us next assume that $\left(B\right)$ holds.
For $r>0$  and $j\in \{1,2\},$  set $B^{j}_{r}= \{ {\bf u} \in \overline{E} : |{\bf u} - {\bf o}_{j}| < r \},$   
$C^{j}_{r} = \{ {\bf u} \in \overline{E} : |{\bf u} - {\bf o}_{j}| = r \},$   
 and let $l^{j}_{r}$ be the length of the image curve $Y(C^{j}_{r});$  
also let $C^{j,\prime}_{r} = G(C^{j}_{r})$  and $B^{j,\prime}_{r}= G(B^{j}_{r}).$  
From the Courant-Lebesgue Lemma, we see that 
for each $\delta\in (0,1)$  and $j\in \{1,2\},$ 
there exists a $\rho_{j}=\rho_{j}(\delta)\in \left(\delta,\sqrt{\delta}\right)$  such that the arclength 
$l_{j,\rho}$  of $Y(C^{j}_{\rho_{j}})$  is less than $p(\delta).$    
\begin{figure}[ht]
\centering
\includegraphics{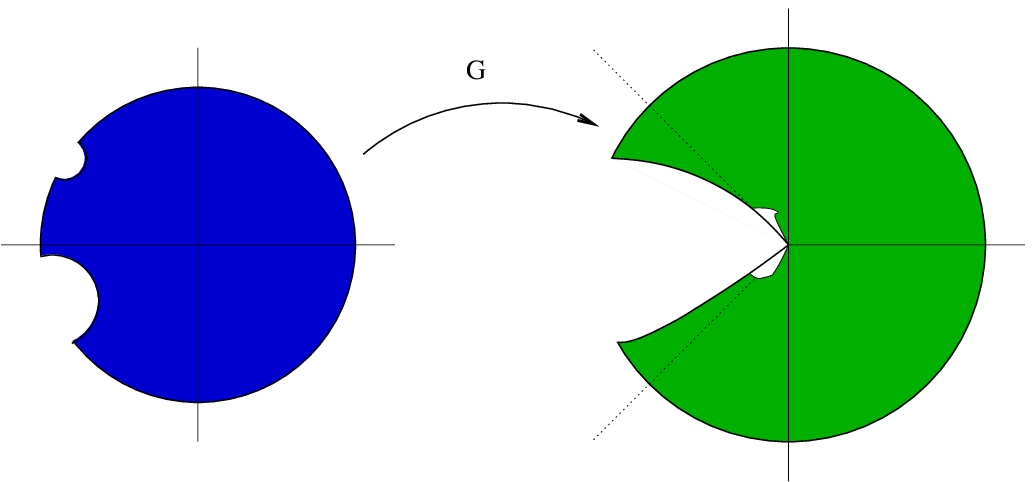}
\caption{$E\setminus \left(\overline{B^{1}_{\rho_{1}(\delta)}} \cup \overline{B^{2}_{\rho_{2}(\delta)}}\right)$  
and $\Omega_{1}$  \label{ONE}}
\end{figure}

We will only consider $\delta\le\delta_{0},$  where $\delta_{0}$  is small enough 
that the endpoints of $C^{j}_{\rho_{j}(\delta)}$  lie on $\sigma_{0}\cup \sigma^{j}_{N}$  for $j\in \{1,2\}$  and 
$C^{1}_{\sqrt{\delta_{0}}} \cap C^{2}_{\sqrt{\delta_{0}}} = \emptyset,$  
where $\sigma^{1}_{N}={\bf o}_{1}{\bf b}$  and $\sigma^{2}_{N}={\bf o}_{2}{\bf a}.$  
For each $\delta\in (0,\delta_{0}),$  the fact that $l_{j, \ \rho_{j}(\delta)}$  is finite for $j\in \{1,2\}$  
implies that  
\[
\lim_{C^{j,\prime}_{\rho_{j}(\delta)}\ni {\bf x}\to {\cal O}}f\left({\bf x}\right) \ \ \ {\rm exists \ for} \ \ j\in \{1,2\}. 
\]
If we set $\Omega_{1}=G\left(E\setminus \left(\overline{B^{1}_{\rho_{1}(\delta)}}\cup \overline{B^{2}_{\rho_{2}(\delta)}} \right)\right)$ 
and define $\phi: \partial\Omega_{1}\to\Real$  by $\phi=f,$   then $\phi$  has (at worst) a jump discontinuity at ${\cal O}.$    
If we consider $\phi$  to be the Dirichlet data for the boundary value problem 
\begin{eqnarray}
\label{PMC2}
{\rm div}(Th) & = &  2H(\cdot,f) \ \ \ \ {\rm in} \ \  \Omega_{1} \\
\label{eq2}
f & = & \phi \ \ \ \  {\rm on} \ \  \partial\Omega_{1}\setminus \{ {\cal O}\},
\end{eqnarray}   
then we may parametrize the graph of $f$  over $\Omega_{1}$  in isothermal coordinates as above and the arguments 
in \cite{EL:86,Lan:88,LS1} can be used to show that $c$  is uniformly continuous on $\Omega_{1}$  and so extends to be continuous 
on $\overline{\Omega_{1}}$   (i.e. Let 
$k:E\setminus \left(\overline{B^{1}_{\rho_{1}(\delta)}}\cup \overline{B^{2}_{\rho_{2}(\delta)}}\right) \to E$  
be a conformal map.  
From \cite{EL:86,Lan:88,LS1}, we see that $c\circ k^{-1}\in C^{0}(\overline{E})$  and so 
$c\in C^{0}\left(\overline{E\setminus \left(B^{1}_{\rho_{1}(\delta)}\cup B^{2}_{\rho_{2}(\delta)}\right) } \right).$)
Since 
\[
\bigcup_{\delta\in (0,1)} \left( E\setminus \left(B^{1}_{\rho_{1}(\delta)}\cup B^{2}_{\rho_{2}(\delta)}\right) \right)=E,
\]
we see $c\in C^{0}\left(\overline{E}\setminus \{{\bf o}_{1}, {\bf o}_{2} \} \right).$

\begin{figure}[ht]
\centering
\includegraphics{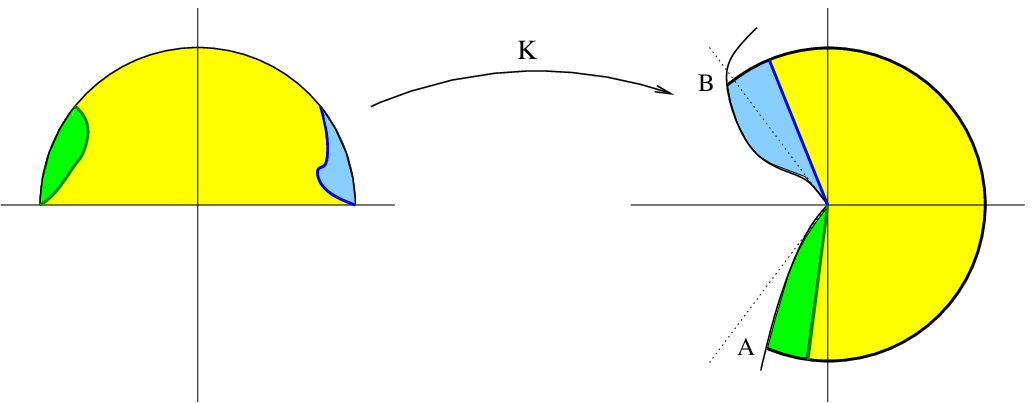}
\caption{ $L(\alpha_{2}),$  $K^{-1}\left( L(\alpha_{2}) \right)$  (blue curves); $L(\alpha_{1}),$  $K^{-1}\left( L(\alpha_{1}) \right)$  (green curves) 
\label{NINE}}
\end{figure}

As at the end of Step 1 of the proof of Theorem 1 of \cite{LS1}, we define $X:B\to\Real^{3}$  by $X=Y\circ g$ and $K:B\to\Real^{2}$  
by $K=G\circ g,$  where $B=\{(u,v)\in\Real^{2} : u^{2}+v^{2}<1, \ v>0\}$  and $g:\overline{B}\to \overline{E}$  is an indirectly 
conformal (or anticonformal) map from $\overline{B}$  onto $\overline{E}$  such that  $g(1,0)= {\bf o}_{1},$   
$g(-1,0)= {\bf o}_{2}$  and $g(u,0)\in {\bf o}_{1}{\bf o}_{2}$  for each $u\in [-1,1].$  
Notice that $K(u,0)={\cal O}$  for $u\in [-1,1]$  (see Figure \ref{NINE}).    
Set $x=a\circ g,$  $y=b\circ g$  and $z=c\circ g,$  so that $X(u,v)=(x(u,v),y(u,v),z(u,v))$  for $(u,v)\in B.$  
Now, from Step 2 of the proof of Theorem 1 of \cite{LS1}, 
\[
X\in C^{0}\left(\overline{B}\setminus\{(\pm 1,0)\}:\Real^{3}\right)\cap C^{1,\iota}\left(B\cup\{(u,0):-1<u<1\}:\Real^{3}\right)
\]
for some $\iota\in (0,1)$  and $X(u,0)=(0,0,z(u,0))$  cannot be constant on any nondegenerate interval in $(-1,1).$  
Define $\Theta(u)= {\rm arg}\left( x_{v}(u,0)+iy_{v}(u,0) \right).$  From equation (12) of \cite{LS1}, we see that 
\[
\alpha_{1}=\lim_{u\downarrow -1} \Theta(u) \ \ \ \ {\rm and} \ \ \ \  \alpha_{2}=\lim_{u\uparrow 1} \Theta(u);
\]
here $\alpha_{1}<\alpha_{2}.$  
As in Steps 2-5 of the proof of Theorem 1 of \cite{LS1}, we see that $Rf(\theta)$  exists when $\theta\in \left(\alpha_{1},\alpha_{2}\right),$ 
\[
\overline{G^{-1}\left( L(\alpha_{2}) \right)} \cap \partial E = \{ {\bf o}_{1} \} \ \   ({\rm and} \ \  
\overline{K^{-1}\left( L(\alpha_{2}) \right)} \cap \partial B = \{ (1,0) \})  \ \  {\rm when} \ \  \alpha_{2}<\alpha
\]
\[
\overline{G^{-1}\left( L(\alpha_{1}) \right)} \cap \partial E = \{ {\bf o}_{2} \} \ \   ({\rm and} \ \  
\overline{K^{-1}\left( L(\alpha_{1}) \right)} \cap \partial B = \{ (-1,0) \})  \ \  {\rm when} \ \  \alpha_{1}>-\alpha
\]
where  $L(\theta)= \{(r\cos(\theta),r\sin(\theta))\in \Omega : 0<r<\delta^{*} \},$
and one of the following cases holds:
\vspace{1mm}

\noindent (a) $Rf$ is strictly increasing or strictly decreasing on $(\alpha_{1}, \alpha_{2})$.  
\vspace{1mm}

\noindent (b) There exist $\alpha_{L}, \alpha_{R}$ so that $\alpha_{1} < \alpha_{L} < \alpha_{R} < \alpha_{2},$  
$\alpha_{R}= \alpha_{L} + \pi$, and $Rf$ is constant on $[ \alpha_{L}, \alpha_{R}]$  and either increasing
on $(\alpha_{1}, \alpha_{L}]$ and decreasing on $[\alpha_{R}, \alpha_{2})$ or decreasing on $(\alpha_{1}, \alpha_{L}]$ 
and increasing on $[\alpha_{R}, \alpha_{2})$.  
\vspace{1mm}

\noindent If $\alpha_{2}=\alpha$  and   $\alpha_{1}=-\alpha,$  then Theorem \ref{CONCLUSION} is proven.  
Otherwise, suppose $\alpha_{2}<\alpha$  and fix $\delta_{0}\in (0,\delta^{*})$  and $\Omega_{0}$  
(see Figure \ref{TWO})  as before in case (i).  
\vspace{2mm}

\noindent {\bf Claim:} Suppose $\alpha_{2}<\alpha.$  Then $f$  is uniformly continuous on $\Omega^{+}_{0},$  where 
\[
\Omega^{+}_{0} \myeq \{(r\cos(\theta),r\sin(\theta))\in \Omega_{0} : 0<r<\delta^{*}, \alpha_{2}<\theta<\pi\}.
\]
\vspace{2mm}

\noindent {\bf Pf:}  Suppose $\alpha-\alpha_{2}<\pi$   (see the blue region in Figure \ref{NINE}).  
Let $\epsilon>0.$  Choose $\delta\in (0,\delta_{0})$  such that  
$p(\delta)+q(p(\delta))<\frac{\epsilon}{4}$  and $p(\delta)<r_{2}-r_{1}.$  
Let $C_{r} = \{ (u,v) \in \overline{B} : |(u,v)-(1,0)| = r \}$  and let $l_{r}$ be the arclength of the image curve $X(C_{r}).$  
The Courant-Lebesgue Lemma implies that for each $\delta\in (0,1),$  there exists a 
$\rho(\delta)\in \left(\delta,\sqrt{\delta}\right)$  such that $l_{\rho(\delta)}< p(\delta).$  
Denote the endpoints of $C_{\rho(\delta)}$  as $(u_{1}(\delta),v_{1}(\delta))$  
and $(u_{2}(\delta),0),$  where $\left(u_{1}(\delta)\right)^{2}+\left(v_{1}(\delta)\right)^{2}=1,$  $v_{1}(\delta)>0$  
and $u_{2}(\delta)\in (-1,1).$
Notice $\Theta\left(u_{2}(\delta)\right)<\alpha_{2};$  let us assume that $\delta$  is small enough that 
$\alpha-\Theta\left(u_{2}(\delta)\right)<\pi.$

\begin{figure}[ht]
\centering
\includegraphics{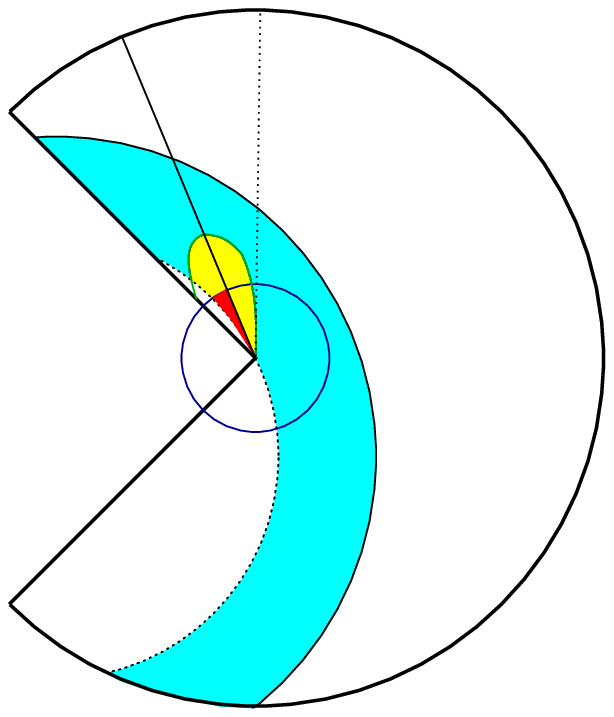}
\caption{$\Omega^{*}\cap {\cal A}\left({\bf p}_{1}\right)$  (blue, yellow \& red regions), 
$\partial B_{\eta(\delta)}({\cal O})$   (blue circle)  \label{FIVE}}
\end{figure}

Now $X\left(C_{\rho(\delta)}\right)$  is a curve whose tangent ray at ${\cal O}$  exists and has direction 
$\theta=\Theta\left(u_{2}(\delta)\right)$  and   
$\partial\Omega\cap X\left(C_{\rho(\delta)}\setminus \{(u_{1}(\delta),v_{1}(\delta)), (u_{2}(\delta),0)\}\right)=\emptyset;$  
hence there exists $\eta=\eta(\delta)>0$  such that 
$\{ {\bf x}\in \Omega^{+}_{0} : |{\bf x}|<\eta(\delta)\}$  (the red region in Figure \ref{FIVE})  is a subset of 
$\Omega_{0} \cap X\left(\{(u,v)\in\overline{B} : |(u,v)-(1,0)|<\rho(\delta)\}\right)$  
(the yellow region plus the red region in Figure \ref{FIVE}).  
From (\ref{Peach}) and the arguments in the proof of the Claim in case (i), we see that $f$  is uniformly continuous 
on $\Omega^{+}_{0}.$

Notice that if $\alpha-\alpha_{2}\ge\pi,$   then we argue as in the proof of the Claim in case (i) and 
see that $f$  is uniformly continuous on $\Omega^{+}_{0}.$  Thus, the Claim is proven.
 
\vspace{2mm}

\noindent Thus $f\in C^{0}\left(\overline{\Omega^{+}_{0}}\right);$  hence (\ref{Tree}) implies  
\[
Rf(\theta)=\lim_{\tau\uparrow \alpha_{2}} Rf(\tau) \ \ \ \ {\rm for \ all} \ \ \theta\in [\alpha_{2},\alpha).
\]  

Suppose $\alpha_{1}>-\alpha.$  Then, as above, $f$  is uniformly continuous on 
\[
\Omega^{-}_{0} \myeq \{(r\cos(\theta),r\sin(\theta))\in \Omega_{0} : 0<r<\delta^{*}, -\pi<\theta<\alpha_{1}\}
\]
and $f\in C^{0}\left(\overline{\Omega^{-}_{0}}\right);$  hence (\ref{Tree}) implies  
\[
Rf(\theta)=\lim_{\tau\downarrow \alpha_{2}} Rf(\tau) \ \ \ \ {\rm for \ all} \ \ \theta\in (-\alpha,\alpha_{1}].
\]
Thus Theorem \ref{CONCLUSION} is proven.

\section{Proof of Theorem \ref{DEMONSTRATION}}  
The parametric representation (\ref{PARAMETRIC}) with properties $\left(a_{1}\right)-\left(a_{5}\right)$  continues 
to be valid and either case $\left(A\right)$  or case $\left(B\right)$  holds true.  

Suppose case $\left(A\right)$  holds.  
Let $q_{1}$  denote the modulus of continuity of the trace of $f$  on the (closed) set 
$\partial^{-}\Omega^{*}$  (i.e. $|f({\bf x}_{1}) - f({\bf x}_{2})| \leq q_{1}(| {\bf x}_{1} - {\bf  x}_{2}|)$  if 
${\bf x}_{1}, {\bf  x}_{2} \in \partial^{-}\Omega^{*}$).  
Fix $\delta_{0}\in (0,\delta^{*})$  with $\sqrt{\delta_{0}}<\min\{|{\bf o}-{\bf a}|, |{\bf o}-{\bf b}|\}.$  
Let ${\bf p}_{1}\in\Real^{2}$  satisfy $|{\bf p}_{1}|=r_{1}$  and $|{\bf p}_{1}-{\bf y}_{1}(\delta_{0})|=r_{1}$  
such that ${\bf p}_{1}$  
lies above (and to the left of) the line through ${\cal O}$  and ${\bf y}_{1}(\delta_{0}).$  
Set $\Omega_{0} = \{ {\bf x}\in \Omega^{*} \ : \ |{\bf x}-{\bf p}_{1}|>r_{1} \}.$
\vspace{1mm}

\noindent {\bf Claim:} $f$  is uniformly continuous on $\Omega_{0}.$
\vspace{1mm}

\noindent {\bf Pf:}   Let $\epsilon>0.$  Choose $\delta\in (0,\delta_{0})$  such that  
$p(\delta)+q(p(\delta))+q_{1}(p(\delta))<\frac{\epsilon}{2}$  and $p(\delta)<r_{2}-r_{1}.$    
Pick a point ${\bf w}\in C^{\prime}_{\rho(\delta)}$  and define $b^{\pm}_{1}:{\cal A}({\bf p}_{1})\to\Real$  by 
\[
b^{\pm}_{1}({\bf x})=f({\bf w})  \pm p(\delta)  \pm h_{{\bf p}_{1}}({\bf x}), \ \ \ \ {\bf x}\in {\cal A}({\bf p}_{1}).
\] 
Notice that 
\[
b^{-}_{1}({\bf x})<f({\bf x})<b^{+}_{1}({\bf x}) \ \ \ \ {\rm for} \ \  
{\bf x}\in B^{\prime}_{\rho(\delta)}  \cap {\cal A}({\bf p}_{1}).
\]
\begin{figure}[ht]
\centering
\includegraphics{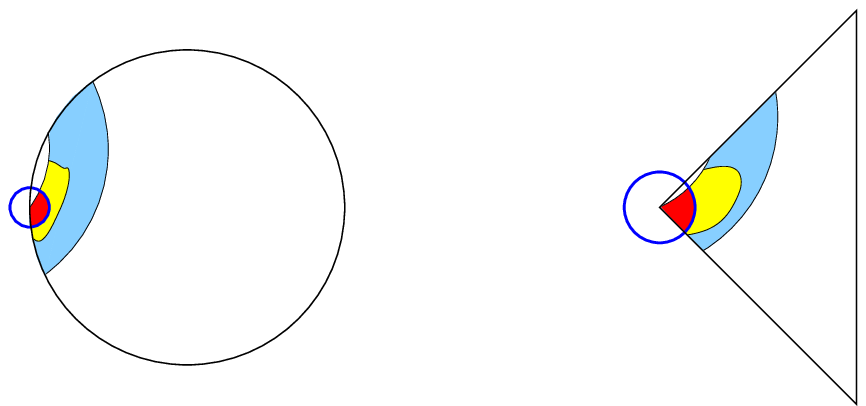}
\caption{$\Omega_{0}\cap {\cal A}\left({\bf p}_{1}\right)$  (blue, yellow \& red regions), 
$\partial B_{\eta(\delta)}({\cal O})$   (blue circle)  \label{SIX}}
\end{figure}

Now there exists $\eta=\eta(\delta)>0$  such that $\{ {\bf x}\in \Omega_{0} : |{\bf x}|<\eta(\delta)\}$  
(the red regions in Figure \ref{SIX})  is a subset of $B^{\prime}_{\rho(\delta)}  \cap {\cal A}({\bf p}_{1})$  
%$\Omega_{0} \cap G\left(\{{\bf w}\in\overline{E} : |{\bf w}-{\bf o}|<\rho(\delta)\}\right)$  
(the yellow regions plus the red regions in Figure \ref{SIX}).
Thus, for ${\bf x}_{1},{\bf x}_{2}\in \Omega_{0}$  satisfying  $|{\bf x}_{1}|<\eta(\delta),$  $|{\bf x}_{2}|<\eta(\delta),$   
we have 
\[
|f({\bf x}_{1})-f({\bf x}_{2})|<2p(\delta)+2q\left(p(\delta)\right)+2q_{1}\left(p(\delta)\right)<\epsilon.
\] 
The remander of the proof of the claim follows as before.
\vspace{1mm}

The proof of Theorem \ref{DEMONSTRATION} in this case now follows the proof of Theorem \ref{CONCLUSION} in the same case.

If case $\left(B\right)$  holds, then the proof of Theorem \ref{DEMONSTRATION} is essentially the same as the proof of 
Theorem \ref{CONCLUSION}; the only significant difference is that $z\in C^{0}\left(\overline{B}\setminus \{(1,0)\}\right)$
(and $c\in C^{0}\left(\overline{E}\setminus \{{\bf o}_{1}\}\right)$)  and hence $Rf(\theta)$  exists for 
$\theta\in [-\alpha,\alpha).$

\end{document}